\documentclass[10pt, reqno]{amsart}
\usepackage{amscd,amsthm,amsfonts,amssymb,amsmath,euscript}
\usepackage[dvips]{graphicx}
\usepackage[dvips]{graphics}
\usepackage[matrix,arrow]{xy}



\newcommand{\subvictor}{\par\medskip \noindent\refstepcounter{subsection}%
{\bf \arabic{section}.\arabic{subsection}.} }

\newcommand{\subsubvictor}{\par\medskip \noindent\refstepcounter{subsubsection}%
{\bf \arabic{section}.\arabic{subsection}.\arabic{subsubsection}.} }

\theoremstyle{definition}

\theoremstyle{remark}

\newcommand{\QQ}{{\mathbb Q}}
\newcommand{\ZZ}{{\mathbb Z}}
\newcommand{\NN}{{\mathbb N}}

\newcommand{\PP}{{\mathbb P}}
\newcommand{\TT}{{\mathbb T}}
\newcommand{\CC}{{\mathbb C}}

\newcommand{\Th}{{\bf Theorem}}
\newcommand{\Lem}{{\bf Lemma}}
\newcommand{\Def}{{\bf Definition}}
\newcommand{\Rem}{{\bf Remark}}

\newcommand{\Corol}{{\bf Corollary}}
\newcommand{\Ex}{{\bf Example}}
\newcommand{\Prop}{{\bf Proposition}}
\newcommand{\Proof}{{\bf Proof}}

\newfont{\smallskob}{cmbx7 scaled\magstep4}
\newfont{\bigskob}{cmbx12 scaled\magstep4}

\newcommand{\itc}[1]{\textup{#1}}

\newcommand{\pic}{\mathrm{Pic}\,}
\newcommand{\spec}{\mathrm{Spec}}

\newcommand{\lcm}{\mathrm{LCM}}

\newcommand{\rd}{\mathrm{det}_\mathrm{right}\,}
\newcommand{\e}{\varepsilon}

\newcommand{\tit}{MINIMAL GROMOV--WITTEN RING}

\begin{document}

\begin{title}
\tit
\end{title}

\author{Victor Przyjalkowski}

\thanks{The work was partially supported by RFFI grant $05-01-00353$ and grant NSh$-9969.2006.1$.}

\address{Steklov Institute of Mathematics, 8 Gubkin street, Moscow 119991, Russia} %




\email{victorprz@mi.ras.ru, victorprz@gmail.com}

\maketitle

\begin{center}
\parbox{340pt}{\small \bf Abstract. \rm
We build the abstract theory of Gromov--Witten invariants of genus 0 for quantum minimal Fano varieties
(a minimal natural (with respect to Gromov--Witten theory) class of varieties).
In particular, we consider ``the minimal Gromov--Witten ring'', i. e. a
commutative algebra with generators and relations of the form used in the Gromov--Witten
theory of Fano variety (of unspecified dimension).
Gromov--Witten theory of any quantum minimal variety is a homomorphism
of this ring to $\CC$.
We prove the Abstract Reconstruction Theorem which states the particular isomorphism of this ring
with a free commutative ring generated by ``prime two-pointed invariants''.
We also find the solutions of the differential equations of type $DN$
for a Fano variety of dimension $N$ in terms of generating series of one-pointed Gromov--Witten invariants.

}
\end{center}

\bigskip
\bigskip
\bigskip

Consider a smooth Fano variety $V$ of dimension $N$. Let $H_H^*(V,\QQ)\subset H^*(V,\QQ)$
be a subspace multiplicatively generated by the anticanonical class $H\in H^2(V,\QQ)$.
It is tautologically closed with respect to the multiplication in cohomology.
The Gromov--Witten theory (more precisely, the set of Gromov--Witten invariants of genus 0)
enables one to ``deform'' the cohomology ring, that is, to define quantum multiplication in the graded space
$QH^*(V)=H^*(V,\QQ)\otimes \CC[q]$. This multiplication
coincides with the multiplication in $H^*(V,\QQ)$ under specification $q=0$.
The subspace $QH_H^*(V)=H^*_H(V)\otimes \CC[q]$ is not closed with respect to the quantum multiplication
in general. The variety is called quantum minimal is the corresponding subspace is closed.

Let $V$ be quantum minimal.
Then there is a natural submodule corresponding to $QH_H^*(V)$ in its quantum $\mathcal D$-module.
The connection in it is given by a matrix of quantum multiplication by the anticanonical
class of $V$ (in the other words, by prime two-pointed genus zero Gromov--Witten invariants).
After regularizing it one get the determinantal operator (the operator of type $DN$) for $V$.
Such operators are deeply studied in~\cite{GS05}.
One of the mirror symmetry conjectures states that this operator coincides with the Picard--Fuchs operator
for dual to $V$ Landau--Ginzburg model (see, for instance~\cite{KM94}).
This is proved, for example, for threefolds,
see~\cite{Pr04} and~\cite{Pr05}.

In the paper we prove the following theorem. Let
$\langle\tau_iH^j\rangle_d$ be the one-pointed Gromov--Witten invariant for the curve of the anticanonical
degree $d$. Consider a regularized $I$-series (the generating series for one-pointed Gromov--Witten invariants)
$$
\widetilde{I}^V=1+\sum_{0\leq j\leq N,d>0} \langle
\tau_{d+j+2}H^{N-j}\rangle_{d} q^{d}h^j/H^N\cdot(h+1)\cdot\ldots\cdot
(h+d)\in \CC[[q]][h]/h^{N}.
$$

\medskip

\Th\ {\bf (\Corol~\ref{corollary:DN's})}. {\it Let
$\widetilde{I}^V=\sum_{0\leq k\leq N}\widetilde{I}^kh^k$, where
$\widetilde{I}^i$ for any $i$
is the series on $q$. Then $N$ functions
$$
\widetilde{I}^0,\ \ \widetilde{I}^0\log (q)+\widetilde{I}^1,\ \
\widetilde{I}^0\log (q)^2/2!+\widetilde{I}^1\log
(q)+\widetilde{I}^2,\ \ \ldots
$$
form the basis of the kernel of the operator of type $DN$ for $V$.
}

\medskip

The idea of the expression of the solutions of equations of type $DN$ in terms
of Gromov--Witten invariants is going back to Witten, Dijkgraaf, and Dubrovin.
The proof of this theorem is based on the relations between Gromov--Witten invariants of genus 0.
This enables us to generalize Gromov--Witten theories of quantum minimal Fanos and their determinantal operators
to the formal theory which does not depend on the variety, and, moreover, on its dimension.
Numerical invariants in this theory become the universal polynomials that are unique for all quantum minimal
Fanos of any dimension.

\medskip

Fix $N\in \NN$. Let $A_N=\CC[a_{ij}]$, $0\leq i,j\leq N$. Let
$\mathcal D=\CC[q,\frac{d}{dq}]$. Put $D=q\frac{d}{dq}$. Consider
the matrix
$$
M=\begin{pmatrix}
  a_{00}(Dq) & a_{01}(Dq)^2 & \ldots & a_{0,N-1}(Dq)^N & a_{0,N} (Dq)^{N+1} \\
  1 & a_{11}(Dq) & \ldots & a_{1,N-1}(Dq)^{N-1} & a_{1,N} (Dq)^N \\
    &   & \ldots &   &   \\
  0 & 0 & \ldots & 1 & a_{NN} (Dq)
\end{pmatrix}
$$
with entries in $A_N\otimes \mathcal D$. Define the operator $L_N\in
A_N\otimes \mathcal D$ by
$$
\rd (D-M)=DL_N \footnote{The operators that come from a symmetric
with respect to the anti-diagonal matrices (i.~e. ones with
$a_{ij}=a_{N-j,N-i}$) are called the operators of type $DN$. }
\label{definition:DN}
$$
(as $A_N\otimes \mathcal D$ is non-commutative, the determinant is
taken with respect to the rightmost column).

\medskip

\Ex. Consider operators $L_N$ as (non-commutative) polynomials in
$q$ and $D$. It is easy to check that the degrees of such operator
in $q$ and $D$ do not depend on its representation as a polynomial.
As a polynomial, $L_2$ (resp.~$L_3$) is of degree $2$ (resp.~$3$) in
$D$ and of degree $3$ (resp.~$4$) in $q$. Write down the analytic
solutions of $L_2\Phi=0$ and $L_3\Phi=0$ in terms of $a_{ij}$'s.

\begin{multline*}
1+a_{00}q+(1/2a_{01}+a_{00}^2)q^2+(7/6a_{01}a_{00}+a_{00}^3+1/3a_{01}a_{11}+2/9a_{02})q^3+ \\
(23/12a_{01}a_{00}^2+1/4a_{01}a_{11}^2+3/8a_{01}^2+1/8a_{22}a_{02}+
5/6a_{11}a_{00}a_{01}+3/16a_{01}a_{12}+  \ \ \ \ \ \ \ \ \ \ \ \ \ \ \ (L_2)\\
43/72a_{02}a_{00}+a_{00}^4+1/6a_{02}a_{11})q^4+O(q^5). \ \ \ \ \  \ \ \ \ \
\end{multline*}

\begin{multline*}
1+a_{00}q+(1/2a_{01}+a_{00}^2)q^2+(7/6a_{01}a_{00}+a_{00}^3+1/3a_{01}a_{11}+2/9a_{02})q^3+\\
(23/12a_{01}a_{00}^2+1/4a_{01}a_{11}^2+3/8a_{01}^2+1/8a_{22}a_{02}+
5/6a_{11}a_{00}a_{01}+3/16a_{01}a_{12}+  \ \ \ \ \ \ \ \ \ \ \ \ \ \ \ (L_3)\\
43/72a_{02}a_{00}+a_{00}^4+1/6a_{02}a_{11}+3/32a_{03})q^4+O(q^5). \ \ \ \ \  \ \ \ \ \
\end{multline*}

\medskip

The first few terms of the solutions coincide.
Moreover, the solutions are ``the same'' if $a_{ij}$'s that are ``out
of bounds'' vanish, i.~e. if we put $a_{i3}=0$ in the solution of $L_3\Phi=0$, we get the solution of $L_2\Phi=0$.

This  same turns out to be true in the general case. Fix two natural
numbers $N_1<N_2$. Let $\Phi_1$ and $\Phi_2$ be the analytic
solutions of $L_{N_1}$ and $L_{N_2}$, normalized by $\Phi_i(0)=1$.
Let $\Phi_2^\prime$ be the series given by putting $a_{ij}=0$ for
$N_1<i,j\leq N_2$ in $\Phi_2$. Then $\Phi_1=\Phi_2^\prime$.
Moreover, it is not difficult to prove (see
Lemma~\ref{lemma:principal minors} and
Proposition~\ref{proposition:Frobenius} below) that for any $n$ if
$N_1\gg n$, then
$$
\Phi_1\mod q^{n}=\Phi_2\mod q^n.
$$

\medskip

A similar statement holds for the logarithmic solutions of $L_N\Phi=0$'s.

So, to get the solutions of such equation, one should ``restrict''
the solutions of the equation of larger index. Moreover, the first
few terms of analitic expansions of different solutions coincide.
Thus, we can define ``the universal series''. ``The restrictions''
of this series to $A_N$ are the solutions of equations $L_N\Phi=0$.
This series is ``the generating series of abstract one-pointed
Gromov--Witten invariants'' in the following sense.

Below we define \emph{the minimal Gromov--Witten ring} $GW$ as a
commutative algebra, with generators and relations of the form used in the Gromov--Witten
theory. Our definition is similar to Dubrovin's definition
of formal Frobenius manifold or Kontsevich--Manin's treatment to the theory of Gromov--Witten invariants.
The difference is that we do not fix
the dimension and consider ``the abstract Gromov--Witten invariants of a Fano variety
of unspecified dimension''. For this we consider ``invariants for the classes with one class replaced by
the Poincar$\mathrm{\acute{e}}$ dual one'' and reformulate Kontsevich--Manin axioms in terms of such ``invariants''.

\medskip

{\bf Definition of the minimal Gromov--Witten ring}.
\label{definition:GW} Consider formal symbols of the form
$$
\langle
\tau_{d_1}H^{i_1},\ldots,\tau_{d_{n-1}}H^{i_{n-1}},\tau_{d_n}H_{r}\rangle,
$$
$n\geq 1$, $i_1,\ldots,i_{n-1},r,d_1,\ldots,d_n\in \ZZ_{\geq 0}$ \itc{(}the last
term is indexed by a subscript!\itc{)}. We write $H^i$, $H_j$ instead of $\tau_0 H^i$, $\tau_0H_j$
for simplicity. Define \emph{the
degree} of such symbol as a number $\sum d_s+\sum i_s-r+(3-n)$. Let $F$ be the set of
all symbols with non-negative degrees. \emph{The
minimal Gromov--Witten ring} is a graded ring
$$
GW=\raisebox{3pt}{$\CC[F]$}\slash\raisebox{-3pt}{$\mathrm{Rel}$},
$$
where $\mathrm{Rel}$ is the ideal generated by the following
relations.
\begin{description}
    \item[GW1 ($S_n$-covariance axiom, cf.~\cite{KM94}, 2.2.1)]
    Consider any permutation $\sigma\in S_{n-1}$. Let $j_k=i_{\sigma(k)}$ and $f_{k}=d_{\sigma(k)}$. Then
$$
\langle
\tau_{d_1}H^{i_1},\ldots,\tau_{d_{n-1}}H^{i_{n-1}},\tau_{d_n}H_{r}\rangle=
\langle
\tau_{f_1}H^{j_1},\ldots,\tau_{f_{n-1}}H^{j_{n-1}},\tau_{d_n}H_{r}\rangle.
$$

    \item[GW2 (normalization, cf.~\cite{KM98}, 1.4.1)] Let $r=\sum d_s+\sum i_s+(3-n)$. Then
$$
\langle
\tau_{d_1}H^{i_1},\ldots,\tau_{d_{n-1}}H^{i_{n-1}},\tau_{d_n}H_{r}\rangle=
\frac{(d_1+\ldots+d_n)!}{d_1!\cdot\ldots\cdot d_n!}\cdot M,
$$
where $M=1$ if $\sum d_j=n-3$ and $M=0$ otherwise.
    \item[GW3 (fundamental class axiom or string equation, cf.~\cite{Ma99}, VI--5.1)]
\begin{multline*}
\langle H^0,\tau_{d_1}H^{i_1}, \ldots,
\tau_{d_{n-1}}H^{i_{n-1}},\tau_{d_n}H_r\rangle=
\sum_{j=1}^{n}\langle
\tau_{d_1}H^{i_1}\ldots\tau_{d_{i-1}}H^{i_{j-1}},
\tau_{d_i-1}H_{i_j}, \tau_{d_{i+1}}H^{i_{j+1}},\ldots, \tau_{d_n}H_r
\rangle,
\end{multline*}
except for the case $\langle H^0,H^i,H_i\rangle$, which is given by {\rm GW2}.

    \item[GW4 (divisor axiom, cf.~\cite{Ma99}, VI--5.4)]
    \begin{multline*}
\langle H^1, \tau_{d_1}H^{i_1},\ldots, \tau_{d_{n-1}}H^{i_{n-1}},
\tau_{d_n}H_r\rangle= d\cdot \langle\tau_{d_1}H^{i_1}, \ldots,
\tau_{d_n}H_r\rangle+ \\
\sum_{s=1}^{n-1} \langle \tau_{d_1}H^{i_1},\ldots,
\tau_{d_s-1}H^{i_s+1},\ldots, \tau_{d_n}H_r\rangle+ \langle
\tau_{d_1}H^{i_1},\ldots, \tau_{d_n-1}H_{r-1}\rangle.
\end{multline*}
where $d>0$ is the degree of the left side.
    \item[GW5 (topological recursion, cf.~\cite{Pa98}, formula 6)]
For any numbers $c_1,\ldots c_n, i_1,\ldots i_n$ and set $S\subset \{1,\ldots, n\}$ denote
the sequence $\tau_{c_{s_1}}H^{i_{s_1}},\ldots,\tau_{c_{s_k}}H^{i_{s_k}}$
($s_1,\ldots,s_k$ are different elements of $S$) by
$\coprod_S$.
For any $n\geq 0$
$$
\langle\coprod_{\{1,\ldots,n\}}, \tau_{d_1}H^{j_1},
\tau_{d_2}H^{j_2}, \tau_{d_3}H_r \rangle= \sum \langle
\tau_{d_1-1}H^{j_1}, \coprod_{S_1}, H_a\rangle \langle H^a,
\coprod_{S_2}, \tau_{d_2}H^{j_2}, \tau_{d_3}H_r\rangle
$$
and
$$
\langle\coprod_{\{1,\ldots,n\}}, \tau_{d_1}H^{j_1},
\tau_{d_2}H^{j_2}, \tau_{d_3}H_r \rangle= \sum \langle
\coprod_{S_1}, \tau_{d_1}H^{j_1}, \tau_{d_2}H^{j_2}, H_a \rangle
\langle H^a, \coprod_{S_2}, \tau_{d_3-1}H_r\rangle,
$$
where the sums are taken over all splittings $S_1\sqcup
S_2=\{1,\ldots,n\}$ and all $a\in \ZZ_{\geq 0}$ such that the degrees of the
symbols in the expression are non-negative \itc{(}notice that the sum is
finite\itc{)}.

\end{description}

\medskip

It turns out that $GW$ has a convenient multiplicative basis.
Consider a ring $A=\CC[a_{ij}]$, $0\leq i\leq j$, $j>0$ and the map $r\colon
A\to GW$ given by $a_{ij}\to \langle H,H^j,H_i\rangle$.

\medskip

\Th~\ref{theorem:reconstruction} {\bf (the Abstract Reconstruction Theorem)}. {\it The map $r$ is an
isomorphism.
}

\medskip

This theorem is an abstract version of the First Reconstruction Theorem of Kontsevich and Manin
(\cite{KM94}, Theorem $3.1$).

Let $a_{00}=0$ (``geometric case'', see Remark~\ref{remark:general case} for the general case).
Via $r$ one may view the coefficients of $L_N\Phi=0$'s and of their solutions as lying in $GW$.

Consider the series
$$
\widetilde{I}=1+\sum_{j\geq 0,i> j-2} \langle \tau_iH_j\rangle
q^{i-j+2}h^j(h+1)\cdot\ldots\cdot (h+i-j+2)\in A\otimes
\CC[[q]][[h]].
$$

In terms of $a_{ij}$'s it rewrites as
\begin{multline*}
\widetilde{I}= 1+\left(
a_{11}h+\left(a_{22}-a_{11}\right)h^2\right)q+ \left(
\frac{a_{01}}{2}+\left( \frac{a_{01}}{4}+\frac{a_{11}^2}{2}+
\frac{a_{12}}{4}\right)h+\right.\\
\left.\left(
-\frac{a_{01}}{8}+\frac{a_{23}}{8}+\frac{a_{11}a_{22}}{2}-
\frac{a_{11}^2}{4}+\frac{a_{22}^2}{4}\right)h^2\right)q^2+O(q^3,h^3).
\end{multline*}

This universal generating series of one-pointed Gromov--Witten invariants defines the solutions of $L_N\Phi=0$
for all $N$. Namely, Theorem~\ref{theorem:reconstruction} and Corollary~\ref{corollary:DN's} directly imply
the following theorem. For any $\CC$-algebra $R$ let $r_N\colon A\otimes R\to A_N\otimes R$
be the map given by $a_{ij}\mapsto a_{ij}$ if $0\leq i,j,\leq N$ and $a_{ij}\mapsto 0$ otherwise.
Define $\widetilde{I}^s\in A\otimes
\CC[[q]]$ by $\widetilde{I}=\sum \widetilde{I}^sh^s$. Define $S_i$'s
by $S_0=\widetilde{I}^0$, $S_1=\widetilde{I}^0\log
(q)+\widetilde{I}^1$, $S_2=\widetilde{I}^0\log
(q)^2/2!+\widetilde{I}^1\log (q)+\widetilde{I}^2$ and so on.

\medskip

\Th~\ref{proposition:virtual}. {\it The set $\{r_N(S_0),\ldots,
r_N(S_{N-1})\}$ is the basis for the space of solutions of the differential equation $L_N\Phi=0$. }

\medskip

Corollary~\ref{corollary:DN's} is the particular case of Theorem~\ref{proposition:virtual}.
The proofs of these two results are almost identical. According to this reason in the first part of the paper
(sections $1$ and $2$) we prove (without using of the abstract Gromov--Witten theory) Corollary~\ref{corollary:DN's}
that have direct applications in the studying of the geometry of Fano varieties.
In the second part (Section $3$) we prove Theorem~\ref{theorem:reconstruction}.
Using this theorem we formally conclude Theorem~\ref{proposition:virtual} from Corollary~\ref{corollary:DN's}.

The paper is organized as follows. In Section $1$ we consider a
quantum minimal Fano variety $V$ of dimension $N$. The two-pointed
Gromov--Witten invariants give the quantum connection and the
differential operator associated with it. The elements of kernel of
this quantum differential operator are given by the $I$-series of
$V$ (i.~e. the generating series of one-pointed Gromov--Witten
invariants of $V$). In Section~$2$ we study regularization of the
quantum differential operator, which gives the operator of type
$DN$. The Frobenius method gives the explicit expressions for the
solutions of $DN$ associated with $V$ in terms of its $I$-series.
All proofs in Section~$1$ and Section~$2$ are based on the
fundamental class axiom, the divisor axiom, and the topological
recursion relations for $V$. In Section $3$ we prove the Abstract
Reconstruction Theorem and the abstract version of theorems of
Section~$1$ and Section $2$, using the same arguments in the
abstract setup. We give the explicit recursive relations for all
solutions of the differential operators in the Appendix.

\medskip

\section{Quantum operators}
\label{section:definitions}

\subvictor {\bf Non-commutative determinants.} Let $R$ be an associative
$\CC$-algebra (not necessary commutative). We consider matrices with
entries in $R$. The indices of matrix $M$ of size $N+1$ run from $0$
to $N$. The submatrix of size $i\times i$ that is the NW corner of $M$ is
called \emph{the $i$-th leading principal submatrix}.

\subsubvictor \Def\ [\cite{GS05}, Definition $1.3$]. {The matrix $M$
with elements in $R$ is called \emph{almost triangular} if
$M_{ij}=0$ for $i+1>j$ and $M_{i+1,i}=-1$. }

\subsubvictor \Def\ [\cite{GS05}, Definition $1.2$].
{\label{definition:right determinant}} {Consider the matrix $M$ with
elements in $R$. The \emph{right determinant} of $M$ is the
determinant taken with respect to the rightmost column:
$$
\rd (M)=\sum_{i=0}^N M_{iN}C_{iN},
$$
where $C_{iN}$ are cofactors taken as right determinants. }

\medskip

For any $(N+1)\times (N+1)$-matrix $M=(M_{ij})_{0\leq i,j\leq N}$ define the matrix $M^\tau$ by
$M^\tau_{ij}=M_{N-j,N-i}$ (i.~e. $M^\tau$ is ``transpose to $M$ with respect to the anti-diagonal'').

\subsubvictor \Lem\ [Golyshev, Stienstra, see~\cite{GS05}, 1.4].
\label{lemma:principal minors} {\it Let $M$ be an almost triangular
$(N+1)\times (N+1)$-matrix. Put
$$
P_0=1,\ \ \ \ \ P_{i+1}=\sum_{j=0}^i M_{ji} P_j.
$$
Then $P_i$ is the right determinant of $i$-th leading principal
submatrix of $M$. In particular, $P_{N+1}=\rd (M)$.
}

\medskip

\Proof. By induction on the size of the submatrices. For $i=1$ this
is trivial. Denote the ($i+1$)-th leading principal submatrix of $M$
by $M_{i+1}$. Notice that the right determinant of the
matrix $M_{i+1}^j$ obtained by deleting last column and $j$-th row
from $M_{i+1}$ equals
$$
\rd M_{i+1}^j =(-1)^{i-j}P_j.
$$
Thus, we have
$$
\rd M_{i+1}=\sum_{j=0}^{i} (-1)^{j+i} M_{ji}\rd
M_{i+1}^j=\sum_{j=0}^i M_{ji} P_j=P_{i+1}.
$$
\qed

\medskip

\subsubvictor \Lem\ [cf. Golyshev, Stienstra, Proposition 1.7
in~\cite{GS05}]. {\it \label{lemma:zero solution} Let $M$ be an almost
triangular matrix and let $\xi=(\xi_0,\ldots,\xi_N)^T$. Let $M\xi=0$.
Then
$$
\rd (M^\tau)\xi_N=0.
$$
}

\medskip

\Proof. We have the following system of equations on $\{\xi_i\}$:
$$
\begin{cases}
    M_{0,0}\xi_0+\ldots+
    +M_{0,N}\xi_{N}=0  \\
    -\xi_0+M_{1,1}\xi_1+\ldots+M_{1,N}\xi_{N}=0  \\
    \ldots \\
    -\xi_{i-1}+M_{i,i}\xi_{i}+\ldots+M_{iN}\xi_N=0 \\
    \ldots \\
    -\xi_{N-1}+M_{NN}\xi_N=0.  \\
  \end{cases}
$$
Let $P_i$ be given by Lemma~\ref{lemma:principal minors} applied to
the matrix $M^\tau$. Let us solve the system moving, step by step,
in the reverse direction. We have $P_{i}\xi_N=\xi_{N-i}$. Thus,
$$
\rd (M^\tau)\xi_N=\left(\sum_{i=0}^{N}
M_{iN}^\tau P_i\right)\xi_N=\sum_{i=0}^{N}M_{0i}\xi_i=0.
$$
\qed

\subvictor {\bf Quantum operators.} Consider a smooth Fano variety
$V$ of dimension $N$ with $\pic (V)\cong \ZZ$. Denote $H=-K_V$ and
$H^*(V)=H^*(V,\QQ)$. (The natural map $\pic (V)\to H^2(V,\ZZ)$ is an
isomorphism for smooth Fanos, so we use the same notation for the
element of $\pic (V)\otimes \QQ$ and for its class in $H^2(V)$.) Let
$H^*_H(V)\subset H^*(V)$ be a divisorial subspace, that is one
generated by the powers of $H$.

Let $\gamma_1,\ldots,\gamma_n\in H^*(V)$ and $d$ be the
anticanonical degree of an effective algebraic curve $\beta\in
H_2(V)$. We denote the respective Gromov--Witten invariant (of genus
zero) with descendants of degrees $d_1,\ldots,d_n\in \ZZ_{\geq 0}$
(see~\cite{Ma99}, VI--2.1)  by $\langle
\tau_{d_1}\gamma_1,\ldots,\tau_{d_n}\gamma_n\rangle_d$.

\subsubvictor The subspace $H_H^*(V)\subset H^*(V)$ is tautologically
closed with respect to the multiplication, i. e. for any
$\gamma_1, \gamma_2\in H^*_H(V)$ the product $\gamma_1\cdot
\gamma_2$ lies in $H_H^*(V)$. The multiplication structure on the
cohomology ring may be deformed.
That is, one can consider \emph{a quantum cohomology ring} $QH^*(V)=H^*(V)\otimes \CC[q]$
(see~\cite{Ma99}, Definition $0.0.2$) with quantum multiplication,
$\star\colon QH^*(V)\times QH^*(V)\rightarrow QH^*(V)$, i.~e. the bilinear map given by
$$
\gamma_1\star \gamma_2=\sum_{\gamma,d} q^d\langle \gamma_1,
\gamma_2, \gamma^\vee\rangle_d \gamma
$$
for all $\gamma_1,\gamma_2,\gamma\in H^*(V)$, where $\gamma^\vee$
is the Poincar$\mathrm{\acute{e}}$ dual class to $\gamma$ (we identify elements of $\gamma\in H^*(V)$
and $\gamma\otimes 1\in QH^*(V)$).
The constant term of $\gamma_1\star \gamma_2$ (with respect to $q$) is
$\gamma_1\cdot \gamma_2$. The subspace $QH^*_H(V)=H^*_H(V)\otimes
\CC[q]$ is not closed with respect to $\star$ in general. The
examples of varieties $V$ with non-closed subspaces $H^*_H(V)$ are Grassmannians $G(k,n)$, $k, n-k>1$ of dimension $>4$
(for instance, $G(2,5)$) or their hyperplane sections of dimension
$\ge 4$.

\subsubvictor \Def. \label{definition:quantum minimality} {The
variety $V$ is called \emph{quantum minimal} if $QH_H^*(V)$ is
quantum closed, i.~e. if for any $\gamma_1, \gamma_2\in H_H^*(V)$,
$\mu\in H_H^*(V)^\bot$ the Gromov--Witten invariant $\langle
\gamma_1, \gamma_2,\mu\rangle_d$ vanishes\footnote{ A Fano variety
is called minimal if its cohomology is as small as it can be (just
$\ZZ$'s in every even dimension). Quantum minimal variety has as
small ``quantum anticanonical part'' as it can be, that is, similar
to the quantum cohomology of minimal one. That's why it is natural
quantum analog of classical minimal variety. }. }

\medskip

In other words, the variety is quantum minimal if and only if
$QH_H^*(V)$ is the subring of $QH^*(V)$.

Throughout the paper we assume $V$ to be quantum minimal.

\subsubvictor Consider a ring $B=\CC[q,q^{-1}]$. Consider the basis
$\{H^i\}$, $i=0,\ldots,N$, of $H_H^*(V)$ (where $H^0$ is the ring
unity and $H=H^1$). Let $H_i$ be the Poincar$\rm \acute{e}$ dual for
$H^i$. Consider the (trivial) vector bundle $HQ$ over $\spec (B)$
with fibers $H_H^*(V)$. Put $h^i=H^i\otimes 1\in H_H^*(V)\otimes B$.
Put $h=h^1$ and $k_V= K_V\otimes 1=-h$. Let $S=H^0(HQ)$. As $S\cong
H^*_H(V)\otimes B$, we can consider quantum multiplication as the
map $\star \colon S\times S\to S$.

Let $D=q\frac{d}{d q}\in \mathcal D=\CC[q,q^{-1},\frac{d}{dq}]$.
Consider a \itc{(}flat\itc{)} connection $\nabla$ on $HQ$ defined on
the sections $h^i$ as
$$
\left( \nabla(h^i),q\frac{d}{dq}\right)=k_V\star h^i
$$
(the pairing is the natural pairing between
differential forms and vector fields).
This connection provides the structure of $\mathcal D$-module for $S$ by
$D(h^i)=(\nabla(h^i),D)$.
Obviously,
\begin{equation*}
\label{equation:quantum differential} D\left(\sum_{i=0}^N
f_i(q)h^i\right) =\sum_{i=0}^N q\frac{
\partial f_i(q)}{\partial q}h^i- h\star \left(\sum_{i=0}^N f_i(q)h^i\right).
\end{equation*}

\subsubvictor \Def. \label{definition:quantum equation} { The
$\CC$-linear operator $D\colon S\rightarrow S$ is called \emph{the
quantum operator}. }

\medskip

Define the operator $D_B\colon S\rightarrow S$ as
$D_B(\sum_{i=0}^N f_i(q)h^i)=\sum_{i=0}^N q\frac{
\partial f_i(q)}{\partial q}h^i$.
Let $h\star h^j=\sum_j \alpha_{ij}h^i$, $\alpha_{ij}\in B$. Define
the matrix $M$ by $M_{ij}=-\alpha_{ij}\in \mathcal
D$ for $i\neq j$ and
$M_{ii}=D-\alpha_{ii}\in \mathcal
D$\footnote{Identify any matrix $A$
with entries in $\mathcal D$ with operator $S\to S$
given by $A(\sum f_i h^i)=\sum_i (\sum_j A_{ij}f_j) h^i$. Then $M$ is the matrix of $D=D_B-h\star$.}.

\subsubvictor \Def. \label{definition:QO} {The differential operator
$L^Q_V=\mathrm{det}_{\mathrm{right}}(M)\in \mathcal D$ is called
\emph{the quantum differential operator} of $V$. . }

\subsubvictor In the following we study flat solutions of
differential equations corresponding to operators we defined. The
solutions are ``formal series with logarithms'' and do not lie in
$S$. So, we need to change the base. Put $T=\CC[[q]][t]/(t^{N+1})$.
Put $B^{form}=B\otimes_{\CC[q]} T$ and $S^{form}=S\otimes_{\CC[q]}
T$. Let $\mathcal D$ act on $'$ via $Dt=1$. So, the informal meaning
of $t$ is $\log (q)$. In the following we consider $D$, $D_B$,
$h\star$, and so on as $\CC$-operators $S^{form}\to S^{form}$ and
$L^Q_V$ as $\CC$-operator $B^{form}\to B^{form}$.

\subvictor {\bf Relations.} For simplicity we use below
Gromov--Witten invariants with negative degree of curve or negative
descendants (which are formally not defined). The convention is that
they equal $0$.

\subsubvictor \Th\ [Topological recursion, see~\cite{Ma99},
VI--6.2.1]. \label{topological recursion}{\it Let
$\gamma_1,\gamma_2,\gamma_3 \in H^*(V)$, $a_1\in \ZZ_{>0}$, $a_2,
a_3, d\in \ZZ_{\geq 0}$. Then

$$
\langle \tau_{a_1}\gamma_1, \tau_{a_2}\gamma_2,
\tau_{a_3}\gamma_3 \rangle_{d}= \sum_{d_1+d_2=d,\ a=0,\ldots,N}
\langle \tau_{a_1-1}\gamma_1, H^a\rangle_{d_1}\langle H_a,
\tau_{a_2}\gamma_2, \tau_{a_3}\gamma_3 \rangle_{d_2}.
$$
}

\subsubvictor \Th\ [The divisor axiom, see~\cite{Ma99}, VI--5.4]. \label{axiom:divisor} {\it Let
$\gamma_1,\ldots, \gamma_n\in H^*(V)$, $\gamma_0=rH\in H^2(V,\QQ)$
be an ample divisor and $a_1, \ldots,a_m \in \ZZ_{\geq 0}$. Then
$$
\langle \gamma_0, \tau_{a_1}(\gamma_1),\ldots,
\tau_{a_m}\gamma_m\rangle_d=
rd\langle\tau_{a_1}\gamma_1, \ldots, \tau_{a_m}\gamma_m\rangle_d+
\sum_{s=1}^m \langle\tau_{a_1}\gamma_1,\ldots,
\tau_{a_s-1}\gamma_0\cdot \gamma_s,\ldots,
\tau_{a_m}\gamma_m\rangle_d.
$$
}

\subsubvictor \Th\ [The fundamental class axiom, see~\cite{Ma99}, VI--5.1]. \label{axiom:fundamental class}{\it
Let $\gamma_1,\ldots,\gamma_k \in H^*(V)$, $a_1, \ldots,a_k \in
\ZZ_{\geq 0}$. Then
$$
\langle\tau_{a_1}\gamma_1, \ldots \tau_{a_k}\gamma_k, H^0
\rangle_d=\sum_{i=1}^k\langle
\tau_{a_1}\gamma_1\ldots\tau_{a_{i-1}}\gamma_{i-1},
\tau_{a_i-1}\gamma_i, \tau_{a_{i+1}}\gamma_{i+1},\ldots,
\tau_{a_k}\gamma_k \rangle_d.
$$
}

\subvictor {\bf Fundamental solution.} Put
$e^{Ht}=\sum_{r=0}^{\infty} \frac{H^rt^r}{r!}\in H_H^*(V)\otimes
B^{form}$ (the sum is finite). Put $\langle
\tau_{d_1}t^{\alpha_1}\gamma_1,\ldots,\tau_{d_s}t^{\alpha_s}\gamma_s\rangle_d=t^{\sum
\alpha_i}\cdot\langle
\tau_{d_1}\gamma_1,\ldots,\tau_{d_s}\gamma_s\rangle_d$. Consider the
matrix $\Phi$ with elements
\begin{multline*}
\Phi_{a}^{b}=\sum_{d\geq 0}q^d\left(\langle \tau_{ d+a-b-1}H^b,
H_a\rangle_d+\langle \tau_{d+a-b}tH^{b+1}, H_a\rangle_d+\langle
\tau_{ d+a-b+1}\frac{t^2}{2}H^{b+2},
H_a\rangle_d+\ldots\right)=\\
\sum_{d\geq 0} q^{d}\langle \tau_{
\bullet}e^{Ht}H^b, H_a\rangle_d,
\end{multline*}
$0\leq a,b\leq N$, where the meaning of $\bullet$ in what follows is the number
$$
N+d-3-\sum_{\mbox{\footnotesize \rm terms}}(\mbox{\rm codimension of
the cohomological class}-1).
$$
We use the notation $\langle \tau_{ \bullet}(H^1+H^2), H_a\rangle_d$
for $\langle \tau_{ d+a-2}H^1, H_a\rangle_d+\langle
\tau_{d+a-3}H^{2}, H_a\rangle_d$ and so on. As two-pointed
Gromov--Witten invariants for the degree zero curve are not defined,
we put $ \langle \tau_{ \bullet}e^{Ht}H^b, H_a\rangle_0=\langle H^0,
\tau_{\bullet}e^{Ht}H^b, H_a\rangle_0. $

\subsubvictor \Prop\ [Pandharipande, after Givental, Proposition 2
in~\cite{Pa98}]. \label{theorem:Pandharipande}  {\it Consider the
sections $\phi^i=\sum_{a=0}^{N} \Phi_{a}^{i} h^a\in S^{form}$
\itc{(}i.~e. those that correspond to the columns of
$\Phi$\footnote{ Informally, $\Phi$ is ``the matrix of fundamental
solutions of equation given by the quantum operator in the standard
basis''.}\itc{)}.

{\bf 1)} The sections $\phi^i$ are flat, i.~e. $D\phi^i=0$.

{\bf 2)}
If $D\phi=0$, then $\phi=\sum_{i=0}^N
\alpha_i\phi^i$, $\alpha_i\in \CC$.
}

\medskip

\Proof\ [Pandharipande].
{\bf 1)}
We need to prove that
$$
D_B\phi^i=h\star \phi^i.
$$
On the left we have
\begin{multline*}
D_B(\sum_{a} \Phi_{a}^{i}h^a)=\sum_{a}\sum_{d\geq 0}(d q^{d} \langle
\tau_\bullet e^{Ht}H^i, H_a\rangle_d+q^{d}\langle \tau_\bullet H\cdot e^{Ht}H^i, H_a\rangle_d)h^a=
\sum_a \sum_{d\geq 0} q^{d}\langle \tau_\bullet e^{Ht}H^i, H, H_a\rangle_d h^a
\end{multline*}
by the divisor axiom~\ref{axiom:divisor}. On the right,
\begin{multline*}
h\star (\sum_a\Phi_{a}^{i}h^a)=\sum_s\sum_{d_1, d_2\geq 0}\sum_a
q^{d_1}\langle \tau_\bullet e^{Ht}H^i, H_a\rangle_{d_1} q^{d_2}
\langle H^a, H, H_s\rangle_{d_2}h^s=
\sum_s \sum_d q^{d\geq 0}\langle \tau_\bullet e^{Ht}H^i, H, H_s\rangle_d
h^s
\end{multline*}
by the topological recursion~\ref{topological recursion}. Both sides
are equal.

{\bf 2)}
The constant term of $\Phi$ (with respect to $t$ and $q$) is the
identity matrix. This means that the columns of $\Phi$ are linearly
independent. The differential operator of order $N+1$
has at most $(N+1)$-dimensional space of solutions, so it is generated by
$N+1$ functions $\phi^i$. \qed

\subsubvictor \Rem. Consider the matrix $M$ (see
Definition~\ref{definition:QO}). The proposition above states that $M\Phi^i=0$
for the column-vectors $\Phi^i=(\Phi_{0}^{i},\ldots,\Phi_{N}^{i})^T$ that
correspond to the sections $\phi^i$.

\medskip

\subsubvictor \Corol. \label{corollary:transpose matrix} {Define the
matrix $\Psi$ by $\Psi_i^j=\sum_{d\geq 0} q^{d}\langle \tau_{
\bullet}e^{Ht}H_i, \ H^{N-j}\rangle_d$, $0\leq i,j\leq N$. Let
$\Psi_i=(\Psi_i^0,\ldots,\Psi_i^N)^T$ be its column-vectors. Then
$M^\tau\Psi_i=0$ for $0\leq i\leq N$. }

\medskip

\Proof. Analogous to the proof of Proposition~\ref{theorem:Pandharipande}.\qed

\subvictor {\bf The solutions.}

\subsubvictor Consider any series $ I=\sum_{i=0}^N I^i(q)h^i\in
\CC[[q]][h]/(h^{N+1})$, $I^0(q),\ldots,I^N(q)\in \CC[[q]]$. Let
$I_r=\sum_{i=0}^r \left(I^{r-i}(q)\frac{t^i}{i!}\right)\in T$.

\medskip

\Def. \label{definition:perturbed solution} {We say that a series
$I$ is \emph{the perturbed solution} of the equation $PI=0$ (or just
the operator $P\in \mathcal D$) if $P I_r=0$ for any $r\leq N$. }

\medskip

In the other words, given $P=P(q, D)$, consider $P_H= P(q,D_B) $
replacing $D$ by $D_B$. Then $I$ is a perturbed solution of $P$ if
and only if $P_H(e^{ht}\cdot I)=0$.

\medskip
Recall that the $I$-series of $V$ is defined by
$$
I^V=1+\sum_{i,j,d\geq 0} \langle\tau_i H_j\rangle_d h^jq^d\in S^{form}.
$$

\subsubvictor \Th. {\it \label{theorem:I-series}

{\bf 1)}
The series $I^V$ is the perturbed solution
of the equation $L^Q_VI=0$.

{\bf 2)} If $L^Q_V I=0$, then
$I=\sum a_i I^V_i$ for some $a_0,\ldots,a_N\in \CC$.

}

\medskip

\Proof. \label{proof:Pandharipande} {\bf 1)} We have
\begin{multline*}
L^Q_V I^V_i=L^Q_V \left(\frac{t^i}{i!}+\sum_{d\geq 0} q^{d}\langle \tau_{
\bullet}e^{Ht}H_i\rangle_d\right)=L^Q_V \left(\sum_{d\geq 0} q^{d}\langle \tau_{
\bullet}e^{Ht}H_i, \ H^0\rangle_d\right)=\\
\rd (M)\Psi_i^N=\rd (M^{\tau\tau})\Psi_i^N=0
\end{multline*}
by Lemma~\ref{lemma:zero solution} and Corollary~\ref{corollary:transpose matrix}.

{\bf 2)} The solutions $I^V_s$ are linearly independent (by
Proposition~\ref{theorem:Pandharipande}), so they form a basis of
$(N+1)$-dimensional space of solutions of differential equation of
order $N+1$ associated with $L^Q_V$. \qed

\section{DN's}
\label{section:regularized}

\subvictor
Let
$$
L^Q_V=P_{V,0}(D)+qP_{V,1}(D)+\ldots+q^nP_{V,n}(D)\in\mathcal D
$$
be a
quantum differential operator of a quantum minimal smooth Fano
variety of dimension $N$ (usually $n=N+1$). Its singularities are
not regular in general.

\subsubvictor \Def\ [see~\cite{Go05}, $1.9$]. {The operator
$$
\widetilde{L}_V=P_{V,0}(D)+qP_{V,1}(D)\cdot
(D+1)+\ldots+q^{n}P_{V,n}(D)\cdot(D+1)\cdot\ldots\cdot (D+n)
$$
is called \emph{the regularization} of $L^Q_V$. }

\medskip

The singularities of all known $\widetilde{L}_V$ are
regular\footnote{They are also regular if the matrix of quantum multiplication by the
anticanonical class is diagonalizable, see~\cite{GS05}, Remark
$3.6$.}. Obviously, $\widetilde{L}_V$ is divisible by $D$ on the
left.

\subsubvictor \Def\ [see~\cite{Go05}, Definition $2.10$]. {The
operator $L_V$ such that $DL_V=\widetilde{L}_V$ is called \emph{the
\itc{(}geometric\itc{)} operator of type $DN$}. }

\medskip

The solutions of equations associated with geometric operators of
type $DN$ are conjectured to be $G$-series\footnote{That is, for any
solution of type $I=\sum a_iq^i$, $a_n\in \QQ$, the following
conditions hold. Let $a_n=\frac{p_n}{q_n}$, $(p_n,q_n)=1$, $q_n\geq
1$. Then $I$ has positive radii of convergence in $\CC$ and
$\overline{\QQ}_p$ for any prime $p$ and there exist a constant
$C<\infty$ such that $\lcm (q_1,q_2,\ldots,q_n)<C^n$ for any $n$.}.

\medskip

Consider any differential operator
$$
P=P_0(D)+qP_1(D)+\ldots+q^nP_n(D)\in \CC[q,q\frac{d}{dq}].
$$
Let
$$
\widetilde{P}=P_0(D)+qP_1(D)\cdot (D+1)+\ldots+q^n
P_n(D)\cdot(D+1)\cdot\ldots\cdot (D+n)
$$
be a regularization as before.

\subvictor {\bf The Frobenius method}. We describe an ``algebraic''
interpretation of the Frobenius method of solving differential
equations. For the standard version see~\cite{CL55}, IV--8.

Let $R=\CC[\e]/(\e^{N+1})$, $N+1\in \NN$. Consider the differential
operator
$$
P_\e=P_0(D+\e)+qP_1(D+\e)+\ldots+q^nP_n(D+\e)\in \mathcal D\otimes
R.
$$

\subsubvictor \Def. {Consider the sequence $\{\overline{c}_i\}$,
$i\geq 0$, $\overline{c}_i\in R$. It is called \emph{a Newton
solution} of $P_\e$, if for any $m\in \ZZ$
$$
\overline{c}_mP_n(m+\e)+\overline{c}_{m+1}P_{n-1}(m+1+\e)+\ldots+\overline{c}_{m+n}P_{0}(m+n+\e)=0
$$
\itc{(}the convention is that $\overline{c}_i$'s with negative subscripts are $0$\itc{)}.

}

\subsubvictor \Prop. \label{proposition:Frobenius} {\it The sequence
$\{\overline{c}_i\}$ is a Newton solution of $P_\e$ if and only if
the series
$$
I=\overline{c}_0+q\overline{c}_1+\ldots\in \CC[[q]]\otimes R
$$
is a perturbed solution of $P$. }

\medskip

\Proof. Recall that $T=\CC [[q]][t]/(t^{N+1})$. We consider $T$ in
the proof as a $\CC$-vector space. Consider the linear space (over
$\CC$)
$$
C=\{(\overline{a}_0,\overline{a}_1,\ldots), \ \overline{a}_i\in R\}
$$
with basis $\{b_{ij}=(0,\ldots,0,\e^{N-j},0,\ldots), i\geq 0, 0\leq
j\leq N\}$ ($\e^{N-j}$ is in the $i$-th place) and the isomorphism
$l\colon C\to T$ given by $b_{ij}\mapsto q^it^j/j!$. It is easy to
see that the formulas $q\cdot b_{ij}=b_{i+1,j}$ and $D\cdot
b_{ij}=(i+\e)b_{ij}$ determine the action of $\mathcal D$ on $C$.
Trivially, $l(q\cdot b_{ij})=q\cdot l(b_{ij})$ and $l(D\cdot
b_{ij})=D\cdot l(b_{ij})$, i.~e. the actions of $\mathcal D$ on $C$
and $T$ commute. Thus, $P(\{\overline{c}_i\})=0$ if and only if
$P(I_N)=0$ (we follow the notations of~\ref{definition:perturbed
solution} with $h=\e$). Analogously, if $P(\{\overline{c}_i\})=0$,
then $P(I_{r})=0$ for $0\leq r\leq N$. \qed

\subsubvictor \Rem. A Newton solution of $P_\e$ exists in $R$ if and only if
$\mathrm{mult}_0 P_0\geq N+1$.

\subsubvictor \Rem. We consider the case $P_0(0)=0$. The cases of
the other roots of $P_0$ are of this type after shifting of variables.

\medskip

\subsubvictor \Corol. \label{theorem:hyperplane principle} {Let
$I=\sum a_{ij}q^i\e^j\in \CC[[q]]\otimes R$ be a perturbed solution
of $P$. Then
$$
\widetilde{I}=\sum a_{ij}q^i\e^j\cdot (\e+1)\cdot\ldots\cdot (\e+i)
$$
is the perturbed solution of $\widetilde{P}$. }

\medskip

\Proof. It follows from the formula for $\widetilde{P}$ and
Proposition~\ref{proposition:Frobenius}.   \qed

\subsubvictor \Corol. \label{corollary:DN's} {Let $V$ be a smooth
quantum minimal Fano variety of dimension $N$ and $L^Q_V$ be the
corresponding quantum differential operator. Let $L_V$ be the
corresponding operator of type $DN$. Define the polynomials $I^i(h)$
in $h$ by $I^V=\sum_{i=0}^\infty I^i(h)q^i$.

{\bf 1)} Let
$$
\widetilde{I}^V=\sum_{i=0}^\infty
I^i(h)\cdot(h+1)\cdot\ldots\cdot(h+i)q^{i}.
$$
Then $\widetilde{I}^V\!\!\!\!\mod h^N$ is the perturbed solution of
$L_VI=0$.

{\bf 2)} If $L_VI=0$, then $I=\sum a_i \widetilde{I}^V_i$ for some
$a_0,\ldots,a_{N-1}\in \CC$. }

\medskip

\Proof. {\bf 1)} By Theorem~\ref{theorem:I-series} $I^V$ is the
perturbed solution of $L^Q_V$. Let $\widetilde{L}_V$ be the
regularization of $L^Q_V$. Then, by
Corollary~\ref{theorem:hyperplane principle}, $\widetilde{I}^V$ is a
perturbed solution of $\widetilde{L}_V$. The relations for the
Newton solution for $L_{V,\e}$ are proportional to the corresponding
relations for $\widetilde{L}_{V,\e}$ modulo $\e^N$ (we identify the
parameter $\e$ in the Frobenius method with $h$). So, the Newton
solutions of $\widetilde{L}_{V,\e}$ and $L_{V,\e}$ coincide modulo
$\e^N$.

{\bf 2)} It follows from
standard arguments on linear independence
(see the proof of theorem~\ref{proof:Pandharipande}). \qed

\subvictor \Ex. The matrix of quantum multiplication for $\PP^ N$ is
$$
\begin{pmatrix}
  0 & 0 & \ldots & 0 & (N+1)^{N+1}q^{N+1} \\
  1 & 0 & \ldots & 0 & 0 \\
    &   & \ldots &   &   \\
  0 & 0 & \ldots & 1 & 0
\end{pmatrix}.
$$
The corresponding quantum differential operator is
$$
L^Q_{\PP^N}=D^{N+1}-(N+1)^{N+1}q^{N+1}.
$$
Let $F$ be a class dual to the hyperplane in $\PP^N$ (so,
$-K_{\PP^N}=(N+1)F$) and $f=F\otimes 1\in S^{form}$. It is easy to
see that the series
$$
I^{\PP^N}=\sum_{d\geq 0}
\frac{q^{(N+1)d}}{(f+1)^{N+1}\cdot\ldots\cdot(f+d)^{N+1}}
$$
is a perturbed solution of $L^Q_{\PP^N}\Phi=0$.

The operator of type $DN$ for $\PP^N$ is
$$
L_{\PP^N}=D^{N}-(N+1)^{N+1}q^{N+1}(D+1)\cdot\ldots\cdot(D+N)
$$
and a perturbed solution of this operator is the series
$$
\widetilde{I}^{\PP^N}=\sum_{d\geq 0}
\frac{q^{(N+1)d}(h+1)\cdot\ldots\cdot(h+(N+1)d)}{(f+1)^{N+1}\cdot\ldots\cdot(f+d)^{N+1}}.
$$

\section{Universality of $DN$'s and solutions}
\label{section:universality}

All formulas above are formal consequences of the formulas
\ref{topological recursion}--\ref{axiom:fundamental class}. So, the
natural idea is to define ``abstract Gromov--Witten theory'', that
is, to consider Gromov--Witten invariants as formal variables with
natural relations. Moreover, if we consider ``invariants'' that
correspond to several classes of type $H^i$ and one
Poincar$\mathrm{}\acute{e}$ dual class of type $H_r$ we may develop
a universal abstract Gromov--Witten theory that do not depend on the
dimension $N$.

A convenient multiplicative basis of $GW$ is given by the First Reconstruction Theorem from~\cite{KM94}.
We follow the notations of the definition of $GW$ on page~\pageref{definition:GW}.

\subvictor \Th\ [{\bf The Abstract Reconstruction Theorem}].
\label{theorem:reconstruction} {\it The map $r\colon \CC[a_{ij}]\to
GW$ is an isomorphism. }

\medskip

\Proof.\ The following relation is the formal implication from the relations of type GW5.

\medskip

\begin{description}
    \item[GW6]
For any finite subset $S\subset \NN$ denote
$H^{i_{s_4}},\ldots,H^{i_{s_k}}$ (where $s_j$'s are distinct
elements of $S$) by $\coprod_S$. Then for any $n\geq 0$
$$
\sum \langle \coprod_{S_1}, H^{i_1},H^{i_2},H_a\rangle \langle
H^a,\coprod_{S_2}, H^{i_3},H_r\rangle= \sum \langle \coprod_{T_1},
H^{i_1},H^{i_3},H_b\rangle \langle H^b,\coprod_{T_2},
H^{i_2},H_r\rangle,
$$
where the sums are taken over all splittings $S_1\sqcup S_2=\{4,\ldots,n\}$,
$T_1\sqcup T_2=\{4,\ldots,n\}$ and all $a$ and $b$ such that the degrees of all symbols are non-negative
\itc{(}notice that both sums are finite\itc{)}.

\end{description}

\medskip

These relations are called quadratic relations in the geometrical case (see.~\cite{KM94}, 3.2.2).
If $S$ is empty, then these relations are called associativity equations or
WDVV equations. Though these relations follows from GW5, we include them in the generators of relations ideal
of $GW$ for simplicity.

\medskip

Let us prove that $r$ is epimorphic.
Our proof is an abstract version of one in~\cite{KM94}. We denote
$r(a_{ij})$ also by $a_{ij}$ for simplicity. We need to prove that any
``invariant'' $\langle \tau_{d_1}H^{i_1},\ldots,\tau_{d_n}H_{r}\rangle$
can be expressed in terms of $a_{ij}$'s.

Applying relation GW4 to one- or two-pointed invariants (the
abstract symbols), we may assume $n\geq 3$ (see the proof of
Proposition $5.2$ in~\cite{Pr04}). Using GW5 (and GW2), we may
assume that all $d_i$'s equal $0$.

Given an invariant $C=\langle H^{i_1+1},\ldots, H^{i_{n}},
H_{r}\rangle$, $n\geq 2$, $i_1>1$, write GW6 for classes $H,
H^{i_1},\ldots, H^{i_{n}}, H_{r}$. We see that, modulo invariants
with lower number of terms, GW4 and GW2, $C$ equals the sum of the
invariants with terms $H^{i_1},\ldots,H^{i_n},H_r$. So, using
GW1--GW6 we may express any invariant in terms of three-pointed
invariants with $H^1$ as the first term, that is, in terms of
$a_{ij}$'s. Thus, $r$ is an epimorphism.

Let us prove that $r$ ia a monomorphism step by step.

{\bf Step 1.} Let GW$3_p$ and GW$4_p$ be the relations of type GW3 and GW4 for invariants without descendants.
Then the ideal $\mathrm{Rel}$ is generated by GW1, GW2, GW$3_p$, GW$4_p$, GW5, GW6 (the relations ``commute'').
Notice that GW$3_p$ is a particular case of GW2.

{\bf Step 2.} Let
$$
GW^\prime=\raisebox{3pt}{$\CC[F^\prime]$}\slash\raisebox{-3pt}{$\mathrm{(GW4_p, GW5, GW6)}$},
$$
where $F^\prime\subset F$ are invariants of positive degree of type
$$
\langle
\tau_{d_1}H^{i_1},\ldots,\tau_{d_{n-1}}H^{i_{n-1}},\tau_{d_n}H_{r}\rangle,
$$
with $i_k\geq i_l$ for $k>l$ and $d_k\geq d_l$ if $i_k=i_l$. (Thus, the left side of any relation of type
GW2 becomes just the notation of the number on the right side.)
Obviously, $GW^\prime\cong GW$.

{\bf Step 3.} Let $A_p=\CC[F^\prime_p]/\mathrm{(GW4_p, GW6)}$, where $F^\prime_p\subset F^\prime$ is the
subset of invariants without descendants.
Let us prove that the natural map $A_p\to GW^\prime$ is a monomorphism.
Consider the order on the invariants, that is, the function $w$ on $F^\prime_p$ given by
$$
w(\langle
\tau_{d_1}H^{i_1},\ldots,\tau_{d_{n-1}}H^{i_{n-1}},\tau_{d_n}H_{r}\rangle)=
(\sum d_j, d_1,i_1,\ldots,d_n,r).
$$
We say that $C_1>C_2$ if $w(C_1)>w(C_2)$ (with respect to the
natural lexicographic order). Define the lexicographic order on the
monomials in $F^\prime$, that is, for any two monomials
$M_1=\alpha\cdot C_1^{a_1}\cdot\ldots\cdot C_n^{a_n}, M_2=\beta\cdot
C_1^{b_1}\cdot\ldots\cdot C_n^{b_n}$ (where $\alpha,\beta\in \CC$
and $C_1>C_2>\ldots>C_n$) say that $M_1>M_2$ if $a_1>b_1$, or
$a_1=b_1$ and $a_2>b_2$, and so on. Denote the leading term of $E\in
\CC[F^\prime]$ with respect to this order by $L(E)$. Denote the
relation of type GW5 with the invariant $C$ on the left side (which
is not uniquely defined!) by $\mathrm{GW5}(C)$. For any prime (i.~e.
without descendants) invariant $C$ put $\mathrm{GW5}(C)=C$. Consider
any $P\neq 0$ in $\CC[F^\prime_p]$ such that $r(P)\in
(\mathrm{GW4_p, GW5, GW6})\lhd \CC[F^\prime]$ for the natural map
$r\colon \CC[F^\prime_p]\to \CC[F^\prime]$. Denote $r(P)$ by $P$ for
simplicity. Let $P=\sum_{j\in J} \beta_j\cdot \prod_{i\in I}
C_i^{b_{i,j}}\mathrm{GW5}(C_j)$ modulo
$(\mathrm{GW4}_p,\mathrm{GW6})$, where $\beta_j\in \CC$. Applying
GW4 we may assume that invariants $C_i$'s contains at least three
terms. That is, we can apply relation of type GW5 to them. Denote
the maximal of the leading terms of all summands of type
$\prod_{i\in I} C_i^{b_{i,j}}\mathrm{GW5}(C_j)$ by $L$. Let
$J_0\subset J$ be the subset of indices such that $L(\prod_{i\in I}
C_i^{b_{i,j}}\mathrm{GW5}(C_j))=L$ for $j\in J_0$. Let $L$ has a
factor with descendants. Then we have
$$
P=\sum_{j\in J_0}\beta_j\cdot \prod_{i\in I} C_i^{b_{i,j}}\mathrm{GW5}(C_j)+(\mbox{summands with smaller leading terms}).
$$
Obviously, $L(\mathrm{GW5}(C))=C$.
One may check that the difference of two relations of type $\mathrm{GW5}(C)$ may be expressed
in terms of relations of type GW5 with smaller leading terms and relations of type GW6.
Thus, expressions of type $\mathrm{GW5}(C_i)$ in the sum on the right side coincides for every $i$ modulo
summands with smaller leading terms.
Then
\begin{multline*}
P=\sum_{j\in J_0}\beta_j\cdot \prod_{i\in I} \mathrm{GW5}(C_i)^{b_{i,j}}\mathrm{GW5}(C_j)+
(\mbox{summands with smaller leading terms})=\\
\sum_{j\in J_0} \beta_j\cdot (\prod_{i\in I\cup J_0}
\mathrm{GW5}(C_i)^{c_i})+(\mbox{summands with smaller leading terms}).
\end{multline*}
As $L(P)<L$, we have $\sum_{j\in J_0} \beta_j=0$.
We get the expression for $P$ with smaller $L$.
Repeating this procedure, we obtain the expression for $P$ with $L(P)=L$, i.~e. without relations of type GW5.
Thus, $P\in (\mathrm{GW}4_p,GW6)$ and $A_p\cong GW^\prime\cong GW$, i.~e.
 invariant may be uniquely expressed in terms of prime ones.

{\bf Step 4.} Let us prove that any prime invariant may be uniquely
expressed in terms of $a_{ij}$'s. We call the invariants of type
$\langle H^k,H^1,\ldots,H^1,H_r\rangle$ trivial since the relations
of type GW6 for them are trivial. Let $F_t\subset F^\prime_p$ be the
subset of trivial invariants. Obviously, $A\cong
A_t=\CC[F_t]/(\mathrm{GW4_p})\cong\CC[F_t]/(\mathrm{GW4_p, GW6})$.
Let $F^\prime_t\subset F^\prime_p$ be the subset of invariants
without terms $H^1$ and GW$6^\prime$ be the relations of type GW6
with invariants with terms $H^1$ replaced by ones without such terms
given by GW$4_p$. Let us prove that $A_t\cong
\CC[F^\prime_t]/(\mathrm{GW6}^\prime)\cong A_p$.

Define the function $w^\prime$ on the elements of $F^\prime_t$ given by
$$
w^\prime(\langle H^{i_1},\ldots, H^{i_{n-1}}, H_{r}\rangle)=(n, i_1,\ldots, i_{n-1},r).
$$
Define the order on monomials in $F^\prime$ and the leading term
$L^\prime(E)$ of any $E\in \CC[F^\prime_p]$ as before. The direct
computation shows that the difference of the two relations of type
GW$6^\prime$ with the same leading terms may be expressed in terms
of the relations of type GW$6^\prime$ with the smaller leading terms
(``the relations of type GW6 commute''). Assume that $P=P(a_{ij})\in
(\mathrm{GW6^\prime})\lhd \CC[F^\prime_t]$. As before, we may obtain
the expression for $P$ in terms of relations of type GW$6^\prime$
containing only trivial invariants. Since these relations vanish,
$P=0$ and $A\cong A_p\cong GW$. \qed

\medskip

\subvictor \Rem. So, the Gromov--Witten theory of a quantum minimal
Fano variety $V$ of dimension $N$ is a particular function from
$GW_N=r(i_N(A_N))$ to $\CC$, where $i_N\colon A_N\to GW$ is given by
$i_N(a_{ij})=a_{ij}$ if $(i,j)\neq (0,0)$ and $i_N(a_{00})=0$.

\medskip

Theorem~\ref{theorem:reconstruction} enables us to define the universal $I$-series $I\in A\otimes \CC[[q]][[h]]$
such that for any $N$ the abstract $I$-series for dimension $N$
$$
I^N=\sum_{i,j} \langle\tau_i H_j\rangle\cdot q^dh^j\in
GW_N\otimes\CC[[q]][h]/h^{N+1}
$$
is the restriction of $I$, that is, $I^N=r_N(I\ \!\!\!\!\mod
h^{N+1})$. Analogously, we may define the universal ``regularized
$I$-series'' $\widetilde{I}$ such that for
$$
\widetilde{I}^N=\sum_{i,j} \langle\tau_i H_j\rangle\cdot
q^dh^j\cdot (h+1)\cdot\ldots\cdot (h+d)\in GW_N\otimes\CC[[q]][h]/h^{N+1}
$$
we have $\widetilde{I}^N=r_N(\widetilde{I}\ \!\!\!\!\mod
h^{N+1})$.

Consider the torus $\TT=\spec\ \CC[q,q^{-1}]$ and the trivial vector
bundle $HQ^N$ with fiber $GW_N\otimes \langle
H^0,H^1,\ldots,H^N\rangle$ ($H^i$'s are just notations for basis
vectors). Let $h^i=1\otimes H^i$. Let
$$
A^N=\begin{pmatrix}
  a_{0,0}q & a_{0,1}q^2 & \ldots & a_{0,N-1}q^N & a_{0,N}q^{N+1} \\
  1 & a_{1,1}q & \ldots & a_{1,N-1}q^{N-1} & a_{1,N}q^N \\
    &   & \ldots &   &   \\
  0 & 0 & \ldots & 1 & a_{N,N}q
\end{pmatrix}
$$
(where $a_{00}=0$).
Define \emph{the abstract quantum connection} $\nabla^N$ by
$$
\left( \nabla^N(h^i),q\frac{d}{dq}\right) =A^Nh^i
$$
(the connection commutes with $a_{ij}$). Repeat all the previous for
the abstract case. In particular, define \emph{the abstract quantum
differential operator} $L^Q_N\in GW_N\otimes \mathcal D$ and the
operator $L_N\in GW_N\otimes \mathcal D$ (recall that after
specialization of abstract Gromov--Witten invariants to the
geometric ones this operator is called geometric $DN$). (It is easy
to see that this operator is the same as one defined on the
page~\pageref{definition:DN}.) Then we obtain the following theorem.

\subvictor \Th. \label{proposition:virtual} {\it

    {\bf 1)} The series $I^N$ is the perturbed solution
    of $L^Q_NI=0$.

    {\bf 2)} The series $\widetilde{I}^N\!\!\!\!\mod h^{N}$
    is the perturbed solution of $L_NI=0$.
}

\medskip

In the other words, given $\widetilde{I}^N$ (resp. $L_N$), one should put $a_{ij}=0$ for $N_0<i,j\leq N$ to obtain
$\widetilde{I}^{N_0}$ (resp. $D^{N-N_0}L_{N_0}$).

\subvictor \Rem. The same holds for $DN$'s. Recall that operator of type $DN$ is $L_N$ with identified
$a_{ij}$ and $a_{N-j,N-i}$. Let $J^N$ be a perturbed solution of $DN$. If $n\ll N$ and $N<N_0$, then
$$
J^N\!\!\!\!\mod (q^n)=J^{N_0}\!\!\!\!\mod (q^n,h^{N}).
$$

\subvictor \Rem.
\label{remark:general case}
Define $L^Q_N$ and $L_N$ as operators in $\CC[a_{ij}]$, $0\leq i\leq j$ (i.~e. let $a_{00}$ be non-zero).
Then the universality for their solutions also holds. The universal series for $L^Q_N$ is $e^{a_{00}q}\cdot I^\prime$,
where $I^\prime$ is given from $I$ by shift $a_{ii}\mapsto a_{ii}-a_{00}$, and the universal series for $L_N$
is the regularization of $e^{a_{00}q}\cdot I^\prime$.

\section{Appendix}

Consider a differential operator $P=\sum_{i=0}^N q^iP_i(D)\in
\mathcal D$. Denote the $r$-th formal derivative of $P$ with respect
to $D$ by $P^{(r)}$.

\subvictor \Th. \label{theorem:relations} {\it The series
$I=\sum_{i=0}^N I^ih^i$ \itc{(}$I^i$'s are series in $q$\itc{)} is a
perturbed solution of $P$ if and only if for any $s\leq N$
$$
\frac{P^{(s)}(I^0)}{s!}+\frac{P^{(s-1)}(I^1)}{(s-1)!}+\ldots+P(I^s)=0.
$$
}

\medskip

\Proof. Notice that
$$
P(tJ(q))=tP(J(q))+P^{(1)}(J(q))
$$
for any $J(q)\in \CC[[q]]$ (see~\cite{BvS95}, Proposition $4.3.1$).
Thus,
$$
P(t^rJ)=\sum_{i=0}^r \binom{i}{r}t^iP^{(r-i)}(J).
$$

For any $s\leq N$
\begin{multline*}
P (I_s)= P
\left(\frac{t^s}{s!}I^0+\frac{t^{s-1}}{(s-1)!}I^1+\ldots+I^s\right)=\sum_{\alpha=0}^s
P\left(\frac{t^\alpha}{\alpha !}I^{s-\alpha}\right)=\\
\sum_{\alpha=0}^s\sum_{\beta=0}^\alpha
\left(\binom{\beta}{\alpha}\cdot
\frac{t^{\beta}P^{(\alpha-\beta)}(I^{s-\alpha})}{\alpha !}\right)=
\sum_{\alpha=0}^s\sum_{\beta=0}^\alpha
\left(\frac{t^{\beta}P^{(\alpha-\beta)}(I^{s-\alpha})}{\beta!(\alpha-\beta)!}\right)=
\sum_{\alpha=0}^s\sum_{\beta=0}^\alpha R_{\alpha,\beta},
\end{multline*}
where
$$
R_{\alpha,\beta}=\frac{t^{\beta}P^{(\alpha-\beta)}(I^{s-\alpha})}{\beta!(\alpha-\beta)!}.
$$

Prove the theorem by induction on $s$. Suppose that it holds for any
$s_0<s$. Then
\begin{multline*}
\frac{P^{(s)}(I^0)}{s!}+\frac{P^{(s-1)}(I^1)}{(s-1)!}+\ldots+P(I^s)=\sum_{a=0}^s\frac{t^a}{a!}
\left(\frac{P^{(s-a)}(I^0)}{(s-a)!}+\ldots+P(I^{s-a})\right)=\\
\sum_{a=0}^{s}\sum_{b=0}^{s-a}
\frac{t^aP^{(b)}(I^{s-a-b})}{a!b!}=\sum_{a=0}^{s}\sum_{b=0}^{s-a}
S_{a,b},
\end{multline*}
where
$$
S_{a,b}=\frac{t^aP^{(b)}(I^{s-a-b})}{a!b!}.
$$
Obviously, $S_{a,b}=R_{a+b,a}$, so
$$
\sum_{a=0}^s\sum_{b=0}^{s-a}S_{a,b}=\sum_{a=0}^s\sum_{b=0}^{s-a}R_{a+b,a}=
\sum_{\alpha=0}^s\sum_{\beta=\alpha}^{s} R_{\beta,\alpha}=
\sum_{0\leq\alpha \leq \beta \leq s} R_{\beta,\alpha}=
\sum_{a=0}^s\sum_{b=0}^{a} R_{a,b}.
$$
Thus,
$$
\frac{P^{(s)}(I^0)}{s!}+\frac{P^{(s-1)}(I^1)}{(s-1)!}+\ldots+P(I^s)=P
(I_s),
$$
which proves the theorem. \qed

\subvictor \Rem. For $s=1$ this theorem is proven in~\cite{BvS95},
Proposition $4.3.2$ and for $s\leq 2$ in~\cite{Tj98}, Appendix B.

\subvictor \Prop\ [{\bf Newton method}].
\label{proposition:prerelations} {\it The series
$$
\Phi=a_0+a_1q+a_2q^2+\ldots\in B,\ \ \ a_i\in \CC.
$$
is a solution of $P\Phi=0$ \itc{(}as a formal series\itc{)} if and
only if for any $m\in \ZZ$
$$
a_mP_N(m)+a_{m+1}P_{N-1}(m+1)+\ldots+a_{m+N}P_{0}(m+N)=0,
$$
where $a_i$'s with negative subscripts are assumed to be $0$.
}

\medskip

\Proof. Straightforward. \qed

\medskip

Theorem~\ref{theorem:relations} and
Proposition~\ref{proposition:prerelations} enable us to find the
relations for the solutions of $P\Phi=0$.

\subvictor \Corol. \label{corollary:recursion relations} {Let
$I=\sum a_{ij}h^iq^j$. Then $I$ is a perturbed solution of $P$ if
and only if for any $s\leq N$ and for any $m\in \NN$
\begin{multline*}
\frac{a_{0,m}P^{(s)}_N(m)+a_{0,m+1}P^{(s)}_{N-1}(m+1)+\ldots+a_{0,m+N}P^{(s)}_{0}(m+N)}{s!}+\\
\frac{a_{1,m}P^{(s-1)}_N(m)+a_{1,m+1}P^{(s-1)}_{N-1}(m+1)+\ldots+a_{1,m+N}P^{(s-1)}_{0}(m+N)}{(s-1)!}+\ldots+\\
a_{s,m}P_N(m)+a_{s,m+1}P_{N-1}(m+1)+\ldots+a_{s,m+N}P_{0}(m+N)=0.
\end{multline*}
}

\bigskip

The author is grateful to V.\,Golyshev for proposing the problem,
explanations and reference on the Frobenius method, to V.\,Lunts and
C.\,Shramov for helpful comments, and to M.\,Kazarian for important
remarks.


\begin{thebibliography}{57}

\bibitem[BvS95]{BvS95}
V.\,Batyrev, D.\,van Straten, {\it Generalized hypergeometric
functions and rational curves on Calabi--Yau complete intersections
in toric varieties}, Comm. Math. Phys. 168, 493--533 (1995).

\bibitem[CL55]{CL55}
E.\,Coddington, N.\,Levinson, {\it Theory of ordinary differential
equations}, McGrow--Hill Book Company inc., New
York--Toronto--London, 1955.

\bibitem[Go05]{Go05}
V.\,V.\,Golyshev, {\it Classification problems and mirror duality},
LMS Lecture Note, ed. N.\,Young, 338 (2007), preprint (2005)
math.AG/0510287.

\bibitem[GS05]{GS05}
V.\,Golyshev, J.\,Stienstra, {\it Fuchsian equations of type $DN$},
preprint (2007) math.AG/0701936.

\bibitem[KM94]{KM94}
M.\,Kontsevich, Yu.\,Manin, {\it Gromov-Witten classes, quantum
cohomology, and enumerative geometry}, Comm. Math. Phys. 164 (1994)
525--562.

\bibitem[KM98]{KM98}
M.\,Kontsevich, Yu.\,Manin, {\it Relations between the correlators
of the topological sigma-model coupled to gravity},
 Commun.Math.Phys. 196 (1998) 385--398.

\bibitem[Ko94]{Ko94}
M.\,Kontsevich, {\it Homological algebra of mirror symmetry}, Proc.
International Congress of Matematicians (Z\"{u}rich 1994),
Birkh\"{a}uzer, Basel, 1995, pp. 120--139.

\bibitem[Ma99]{Ma99}
Yu.\,Manin, {\it Frobenius manifolds, quantum cohomology, and moduli
spaces}, Colloquium Publications. American Mathematical Society
(AMS). 47. Providence, RI: American Mathematical Society (AMS)
(1999).

\bibitem[Pa98]{Pa98}
R.\,Pandharipande, {\it Rational curves on hypersurfaces [after
A.\,Givental]}, Societe Mathematique de France, Asterisque. 252,
307--340, Exp. No.848 (1998).

\bibitem[Pr04]{Pr04}
V.\,Przyjalkowski, {\it Gromov--Witten invariants of Fano threefolds
of genera $6$ and $8$}, preprint (2004), math.AG/0410327.

\bibitem[Pr05]{Pr05}
V.\,Przyjalkowski, {\it Quantum cohomology of smooth complete
intersections in weighted projective spaces and singular toric
varieties}, preprint (2005), math.AG/0507232.

\bibitem[Tj98]{Tj98}
E.\,Tjotta, {\it Rational curves on the space of determinantal nets
of conics}, Doctoral Thesis, preprint (1998) math.AG/9802037.


\end{thebibliography}
\end{document}